\documentclass[reqno]{amsart}

\usepackage{amscd}
\usepackage{amssymb}
\usepackage{tikz}

\theoremstyle{plain}
\newtheorem{theorem}{Theorem}[section]
\newtheorem{lemma}[theorem]{Lemma}

\theoremstyle{definition}

\newtheorem{def-prop}[theorem]{Definition-Proposition}
\newtheorem{example}[theorem]{Example}

\theoremstyle{remark}
\newtheorem{remark}{Remark}[section]
\newcommand{\Chi}{\mathcal X}
\newcommand{\Rl}{\mathbb R}
\begin{document}
\title[MLE Existence Algorithm]{Algorithm to check Maximum Likelihood Estimate Existence for integrated PCA}
\author{D.A.Shmelkin}
\email{dmitri.shmelkin@gmail.com}
\begin{abstract}
Being encouraged by \cite{am} that provides an amazing bridge between Statistics and Invariant Theory, and especially by \cite{fm}, where quiver semi-invariant techniques apply to verify the existence of MLE for a recent iPCA model, we provide an enhancement to \cite{fm}. Our  Theorem \ref{main_result}  yields necessary and sufficient conditions for MLE to exist generically for any dimension vector. The conditions can be easily checked with either of our software \cite{t} or \cite{ff} based on Derksen-Weyman algorithm and simplifying the application for statistics practitioners and non-specialists in quivers. For those deep in quiver Representation Theory, Theorem \ref{main_result} relates the MLE existence to the local semi-simplicity of representations as introduced in \cite{sh07}. Moreover, we implemented in \cite{ff} the Flip-Flop algorithm from \cite{ta} and checked that it does work in all cases guaranteed by Theorem \ref{main_result}. We hope that our elementary and short text can serve for the experts in both domains as a warm start in a new category. And we welcome practitioners to use our code at \cite{ff} for their own experiments.
\end{abstract}
\maketitle
\tableofcontents

\section{Gaussian Group Model}
In their amazing work \cite{am} Amendola et. al. established a strong connection between Statistics and Algebraic Transformation Groups as follows. Observing a random vector $\Chi\in \Rl^p$ as a set of $n$ samples constituting a $p\times n$ matrix $X$, a {\it Gaussian multivariate model} is often proposed such that  $X$ is thought of as sampled from a distribution $\Chi = g {\mathcal I}$, where ${\mathcal I}\in \Rl^p$ is a $p$-tuple of standard independent Gaussian variables with mean 0 and variance 1 and $g\in GL(p,\Rl)$ is some hidden invertible transformation of $\Rl ^p$.  
In this situation a self-adjoint positive definite operator $gg^\top$ is the covariance of $\Chi$, ${\mathcal F}=(gg^\top)^{-1}$ is the concentration, and the model applies as an estimator of ${\mathcal F}$ given the observed $X$, $F=(\frac{1}{n}\sum_{j=1}^n X_j X_j^{\top})^{-1}$, where $X_j$ is the column of $X$ corresponding to $j$-th observation of $\Chi$. Notice that each summand, $X_j X_j^{\top}$ is a rank 1 $p\times p$ matrix and the sum can't be invertible if $n<p$. But if $n\geq p$ and the sum is invertible, then the above formula for $F$ is a Maximum Likelihood Estimate of ${\mathcal F}$.

The authors of \cite{am} went further and proposed to consider a {\it Gaussian Group Model} $\Chi = g {\mathcal I}$, such that the operator $g$ a priori belongs to some algebraic subgroup $G\subseteq GL_p({\mathbb K})$ where either ${\mathbb K}=\Rl$ or ${\mathbb K}={\mathbb C}$ and ${\mathbb K}^p$ is equipped with a positive defined scalar product, $\Vert v\Vert ^2, v\in {\mathbb K}^p$, Hermitian in the case ${\mathbb K}={\mathbb C}$. Without loss of generality we may assume $\Vert v\Vert ^2= {\bar v}^{\top}v$. It is additionally important to require $G$ to contain all ${\mathbb R}_{>0}$-scalar linear operators in $GL_p({\mathbb K})$ and $G$ self-adjoint with respect to the scalar product, i.e., the adjoint operator  $g^*= {\bar g}^{\top}\in G$ for $g\in G$.

As a known example of such situation, they mention {\it Matrix Normal Model} that applies to the case, where observations are in fact matrices not just vectors. More precisely, $p=p_1\cdot p_2$ so that $\Rl^p$ can be thought of as a tensor product $\Rl^{p_1}\otimes \Rl^{p_2}$ and the the group $G$ is $GL(p_1,\Rl)\times GL(p_2,\Rl)$ naturally embedded in $GL(p,\Rl)$. For this model, the existence of MLE was explored and algorithm proposed in \cite{d}.

\section{Existence of MLE, closed orbits, and Kempf-Ness criterion}
Let $V={\mathbb K}^p$ be a vector space and $G\subseteq GL(V)$ be an algebraic group. For $n$ observed vectors, constituting the $p\times  n$ matrix $X$ we are looking for a Maximum Likelihood Estimate $g\cdot X, g\in G$.
From the equivariant point of view, $X$ belongs to the space $V^n=V\oplus \cdots \oplus V$ of $n$-tuples of vectors in $V$ and the norm naturally extends to $V^n$ so that we have: 
 
\begin{equation}\label{norm_tuple}
\Vert g\cdot X\Vert = \sum_{j=1}^n \overline{(gX_j)}^{\top}gX_j=\sum_{j=1}^n \overline{X_j}^{\top} g^*g X_j.
\end{equation} 

The main observation of \cite{am} relates the likelihood estimation to the norm minimization over the orbit. Let $\widetilde{SL}(V)\subseteq GL(V)$ be the subgroup of operators with determinant $\pm 1$. Assume that $G$ is self-adjoint and contains ${\mathbb R}_{>0}$ scalar operators:

\begin{lemma}\label{main_am} \cite[Proposition~3.4]{am} Let $X$ be a $p\times n$ matrix of samples. If  $G$ contains ${\mathbb R}_{>0} \cdot Id$, then the MLE, if exists, has the form $\lambda h^*h$ where:

{\bf 1.}  the norm $m=\Vert h\cdot X\Vert$ is minimal over all $h\in G\cap \widetilde{SL}(V)$.

{\bf 2.} $\lambda$  maximizes the function $\lambda m/n - p\log \lambda$ over all $\lambda\in \Rl_{> 0}$.
\end{lemma}
In particular, the norm $\Vert h\cdot X\Vert$ reaches minimum over the orbit $(G\cap \widetilde{SL}(V))X$ and, by Kempf-Ness Criterion, \cite{kn}, the  $G\cap SL(V)$-orbit of $X$ is Zarisski  closed in $V^n$ or, in the terminology from \cite{am}, $X$ is {\it polystable}. Moreover, by \cite[Theorem~3.10]{am} MLE for $X$ exists if and only if the $G\cap SL(V)$-orbit of $X$ is closed; if additionally the stabilizer is a finite group (such $X$ are called {\it stable}), then MLE is unique.

\begin{example} Consider the standard PCA as an application of Lemma \ref{main_am}. Here $G$ is the whole of $GL(V)$. It is well known that the action of $SL(V)$ on $V^n$ has a closed orbit at $X\in V^n$ if and only if the rank of $X$ is $p$; in this case the stabilizer of $X$ consists of just identity. These results are inline with the Gaussian multivariate model above: MLE exists and is unique if and only if ${\rm rk} X=p$. So assume $n\geq p$, and, according to  Lemma \ref{main_am}, the MLE is a multiple of $Q=h^*h$ such that $\det(h)=\pm 1$ and the length of $h\cdot X$  takes minimum. Applying the well-known property of trace to (\ref{norm_tuple}), we get:
\begin{equation}\label{norm_as_trace}
\Vert h\cdot X\Vert = \sum_{j=1}^n \overline{X_j}^{\top} Q X_j= \sum_{j=1}^n {\rm Tr}X_j  \overline{X_j}^{\top}Q= {\rm Tr} SQ, 
\end{equation}
where $S= \sum_{j=1}^n X_j  \overline{X_j}^{\top}.$ For a positive defined symmetric $p\times p$ matrix, $R$, all eigenvalues are positive. Hence, the inequality for the arithmetic and geometric mean of the eigenvalues can be written as ${\rm Tr}R\geq p  \det R^{1/p}$. In particular, since $\det Q=1$, then ${\rm Tr} QS={\rm Tr} Q^{1/2} S Q^{1/2}\geq p \det  S^{1/p}$. And the matrix $Q_0=\det S^{1/p} S^{-1}$ has $\det Q_0=1$ and ${\rm Tr} Q_0S= p \det  S^{1/p}$. So we retrieve $S^{-1}$ as a multiple of MLE. \qed
\end{example}

\section{Quivers representations}
Here we need to briefly recall the concept of {\it quiver} $Q$ as a directed graph with the set $Q^0$ of vertices and $Q^1$ of arrows such that an arrow $\varphi\in Q^1$ goes from the tail $t\varphi\in Q^0$ to the head $h\varphi\in Q^0$. A representation $V$ of $Q$ is an assignment of a vector space $V(a)$ to every vertex and a linear map $V(\varphi),V(t\varphi)\to V(h\varphi)$, to every arrow, the dimension $\dim V=\dim V(a),a\in Q^0$.
For a dimension  vector $\alpha\in {\mathbb Z}_{\geq 0}^{Q^0}$ let  $R(Q,\alpha)=\oplus_{\varphi\in Q^1} {\rm Hom} ({\mathbb K}^{\alpha[t\varphi]},{\mathbb K}^{\alpha[h\varphi]})$ be the vector space of all representations of $Q$ of dimension $\alpha$. Every representation allows for a decomposition in a sum of indecomposable summands, unique up to isomorphisms and permutations of factors. The dimension vectors of indecomposables constitute, as the famous Kac Theorem says, \cite{kac}, a (positive) root system corresponding to the underlying undirected graph of $Q$. In particular, a root $\beta$ is either {\it real} such that $\langle \beta,\beta\rangle=1$ or {\it imaginary} such that $\langle \beta,\beta\rangle\leq 0$, where 
\begin{equation}\label{euler}
\langle \alpha,\beta\rangle= \sum_{i\in Q^0} \alpha[i]\beta[i] - \sum_{\varphi\in Q^1} \alpha[t\varphi]\beta[h\varphi]
\end{equation}
is the Euler bilinear form naturally related to $Q$. Kac also introduced the concept of {\it generic decomposition}\footnote{Kac initially called the decomposition {\it canonical} but after several authors we prefer the word generic as better reflecting the meaning and non-overlapping with another canonical decomposition for tame quivers.} of dimension vector into the sum of roots corresponding to the indecomposable summands for a generic representation.  Following A.Schofield ideas, H.Derksen and J.Weyman proposed in \cite{dw} a quite efficient algorithm to compute it for the case, where $Q$ is a DAG. We implemented that algorithm and several others in \cite{t} and published at \cite{sh09}. \cite{t} is a standalone software with graphical user interface that allows doing computations with quiver representations.

The group $GL(\alpha)=\prod_{j\in Q^0,\alpha[j]>0}GL_{\alpha[j]}$ acts naturally on $R(Q,\alpha)$ and the orbits are the isomorphism classes of representations. If $Q$ is a DAG, then regular invariant functions, ${\mathbb K}[R(Q,\alpha)]^{GL(\alpha)}$ are just constant but rational invariant functions exist and the field ${\mathbb K}(R(Q,\alpha))^{GL(\alpha)}$ is generated by ratios $f/g$, where both $f,g$ are $GL(\alpha)$-semi-invariants of the same weight, $\chi_{\sigma}(g)=\prod\det(g_i)^{\sigma[i]}$, where $g=(g_i,i\in Q^0)\in GL(\alpha)$, $\sigma\in {\mathbb Z}^{Q^0}$. 

\section{Matrix Normal Model from Quiver perspective} Based on Lemma \ref{main_am},
the authors of \cite{dm} considered the Matrix Normal Model from the perspective of Invariant Theory for representations of  Kronecker quiver with two vertices, $n$ arrows of the same directions, and of dimension $\alpha=(p_1,p_2)$.  Namely in \cite{dm}, MLE for Matrix Normal Model is investigated and, for any $p_1,p_2$, precise minimal values of $n$ are found such that a generic $X\in \Rl^{p_1}\otimes \Rl^{p_2}$ is polystable or even stable, see \cite[Theorem~1.2]{dm}, or in other words, MLE exists or exists and unique, respectively. For example, if $n\geq \frac{p_1}{p_2}+\frac{p_2}{p_1}$, then MLE exists. If we compare this with the requirement $n\geq p_1p_2$ for MLE in PCA to exists, the Matrix Normal Model has therefore more applications because requires less samples for MLE to exist, for example $n=3$ is enough if the ratio $p_1/p_2\in((3-\sqrt{5})/2,(3+\sqrt{5})/2)$. This advantage justifies using a more involved model.

\section{Integrated PCA model}
Independently of \cite{am}, another interesting model, iPCA, was proposed in \cite{ta} such that the observations constitute  a $n\times p$ matrix $X$, similar to PCA (but in a transposed form) such that each of $n$ rows is therefore a vector of $p$ features observed but the features are grouped, with group sizes $p_1,p_2,\cdots,p_k, k\geq 1$ such that $p=p_1+p_2+\cdots+p_k$. Let $X_j$ denote the $n\times p_j$ submatrix of $X$ corresponding to that group of $p_j$ features, $j=1,\cdots,k$, so we can write $X=(X_1,X_2\cdots,X_k)$.
In difference to PCA (and to the Matrix Normal Model), the {\it Integrated PCA} model explores the diversity not only among features, but among samples as well. In particular, the iPCA model for $k=1, p_1=p$ is not the same as standard PCA that only explores the diversity of features as represented by samples. Here, the iPCA model looks for a $n\times n$ invertible matrix $\Sigma$ and $p_j\times p_j$ invertible matrices $\Delta_j,j=1,\cdots,k$ such that for every $j=1,\cdots,k$ the {\it whithened} data
\begin{equation}
Z_j=\Delta_j^{-1/2}X_j\Sigma^{-1/2}
\end{equation}
can be regarded as sampled from standard Gaussian $p_j\times n$ matrix distrubution (see \cite[formula~(5)]{ta}). The authors of \cite{ta} gave the definition of the model but their conclusion about the existence of MLE is not very optimistic. Namely, \cite[Theorem~5]{ta} claims that for generic $X$ (${\rm rank} X = n, {\rm rank} X_j=p_j,j=1,\cdots,k$) MLE does not exist unless $p_j=n$ for every $j=1,\cdots,k$. Below we prove that actually MLE exists in many cases not covered by the condition from \cite{ta} so there should be some explanation for this in \cite{ta} that we weren't able to find so far. The rationale of \cite{ta} is that for MLE to exist, a regularization is necessary and they define a {\it penalized MLE}, prove existence and yield a Flip-Flop algorithm. Amazingly, they also give a very simple Flip-Flop Algorithm for usual MLE, \cite[Algorithm~2]{ta} that works once MLE exist. As we see below, MLE exists in many cases and the algorithm from \cite{ta} indeed converges for all such cases we tested! In Section \ref{FFS} we show this simple algorithm and discuss our experiments with it,

Using the ideas from \cite{am} in the version from \cite{dm}, the authors of \cite{fm} related the MLE for iPCA to the stability of representations of star quiver $Q_k$ with $k+1$ vertices $[0,1,\cdots,k]$ and $k$ arrows $j\to 0$.  Set $\alpha=[n,p_1,\cdots,p_k]$ and consider the group $GL(\alpha)$ acting on  $R(Q_k,\alpha)$. Regarding $X$  as a point in $R(Q_k,\alpha)$, the action is as follows:
\begin{equation}
(g,h_1,\cdots,h_k)X= g X diag(h_1,\cdots,h_k)^{-1}.
\end{equation}
As one may observe in Lemma \ref{main_am}, $n$, the number of samples, plays the role of the number of copies of a vector space, $V$, acted upon by $G$. As we mentioned above, samples in iPCA aren't treated independently as in PCA. That is why, Lemma \ref{main_am} applies to $V=R(Q_k,\alpha)$ and the number of copies equal to just 1! At first glance, it sounds  like in iPCA we have just 1 sample and it may seem not  statistically significant. But this "single" sample allows to address the statistical inference problem, to which the Lemma applies. As a take away of \cite{fm} we give a statement that we slightly reformulate in our terms:

\begin{lemma}\label{fm-equiv} \cite[2.6]{fm}
Let $\sigma=[-p,n,n,\cdots,n]$ be the weight of a character $\chi_{\sigma}=\det(g)^{-p}\prod_{j=1,\cdots,k}\det{h_j}^n$ of $GL(\alpha)$. Let $GL(\alpha)_{\sigma}$ be the kernel of $\chi_{\sigma}$. Then an MLE for $X$ exists if and only if $X\neq 0$ and $GL(\alpha)_{\sigma}$ orbit of $X$ is closed in $R(Q_k,\alpha)$. Moreover, if additionally the stabilizer $(GL(\alpha)_{\sigma})_X$ is finite, then MLE is unique.   Conversely, if MLE is unique, then $(GL(\alpha)_{\sigma})_X$ is compact.
\end{lemma}

Based on Lemma \ref{fm-equiv}, a set of results regading existence and uniqueness of MLE for {\it generic} $X\in R(Q_k,\alpha)$ are provided in \cite{fm}, and the most explicit ones are \cite[Theorem~1.1]{fm} that yields a {\it sufficient} condition for MLE to generically exist uniquely, while in \cite[Theorem~1.2]{fm} the case of $p_1=p_2=\cdots p_k=q$ is considered and necessary and sufficient condition for MLE to generically exist or exist uniquely are given for any $(n,q,k)$. However, the general case of unequal $p_j$ is not covered by explicit criteria or simply programmable algorithms. 

The quiver representations with a closed $GL(\alpha)_{\sigma}$-orbit for a character $\sigma$ are called {\it $\sigma$-polystable} in \cite{am},\cite{dm},\cite{fm} and this has a clear relation to the theory of $\sigma$-stable or polystable quiver representation that goes back to the brilliant A.D.King paper \cite{k}.  In \cite{sh07} we called such representations
{\it locally semi-simple} and provided  equivalent descriptions of such representations.
Using that definition, we claim:

\begin{theorem}\label{main_result}
For  $\alpha=[n,p_1,\cdots,p_k]$ set $p=p_1+\cdots+p_k$. Let $\alpha=m_1\beta_1+\cdots+m_l\beta_l$ be the generic decomposition such that $\beta_i\neq \beta_j$ for $i\neq j$. The following conditions are equivalent:

{\bf a.} MLE exists for a generic $X\in R(Q_k,\alpha)$

{\bf b.} Generic $X\in R(Q_k,\alpha)$ is locally semi-simple

{\bf c.} For every $j=1,\cdots,l$ holds $p\beta_j[0]=n(\beta_j[1]+\cdots+\beta_j[k])$.

{\bf d.}  $\langle \beta_i,\beta_j\rangle=0$ for $i\neq j$.
\end{theorem}

This theorem yields a practical efficient method to verify, whether or not, iPCA problem for a
generic measurement $p\times n$ matrix $X$ with feature groups of sizes $p_1,\cdots,p_k$ allows an MLE. Indeed, according to it, it is enough to compute the generic decomposition and check the condition {\bf c}. And to compute generic decomposition, one may use the wonderful Derksen-Weyman algorithm from \cite{dw}. Moreover, in \cite{t} we gave an
implementation for Derksen-Weyman algorithm and in  \cite{sh09} enhanced it to a one for generic locally semi-simple representation also implemented in \cite{t}.

\begin{remark} The equivalence of conditions {\bf a} and {\bf c} follows from Theorem 4.3 and Lemma 5.2 in \cite{fm} but in that paper the authors didn't mention the Derksen-Weyman algorithm  as a method to convert this criterion into a practical tool available to statistics practitioners and non-specialists in Quiver Representations.
\end{remark}
\begin{remark}
The dimension vectors $\alpha$ fulfilling the equivalent conditions of Theorem \ref{main_result} exist and not only the Schurian roots such that the generic decomposition is trivial. To illustrate it as an example of \ref{main_result} application, take $k=5$, i.e., consider iPCA with 5 groups of features. In terms of quiver, this is  $Q_5$, a {\it wild} quiver. Take $\alpha=[5,4,3,1,1,1]$, so $n=5,p=10$. Applying Derksen-Weyman algorithm \cite{dw} or running either of our software  \cite{t},\cite{ff} we get the generic decomposition as follows:
\begin{equation}
[5,4,3,1,1,1]=[3,2,1,1,1,1]+2\times[1,1,1,0,0,0]
\end{equation}
Then we check the condition {\bf c}: $10\times 3= 5(2+1+1+1+1)$ and  $10\times 1= 5(1+1+0+0+0)$ hold true both. Hence, a generic  matrix $X$ of just 5 observations in ${\mathbb K}^{10}$, with groups of features of sizes 4,3,1,1,1 allows an MLE in iPCA model.
\end{remark}

\section{Proof of Theorem \ref{main_result}}
 
Let $X\in R(Q_k,\alpha)$ be generic. By Lemma \ref{fm-equiv}, If MLE exists for $X$, then $GL(\alpha)_{\sigma}$-orbit of $X$ is closed, so {\bf a} implies {\bf b}.  That {\bf a} is equivalent to {\bf c} is proven in Theorem 4.3 and Lemma 5.1 and 5.2 in \cite{fm}. Now assume {\bf b} and let $X=m_1V_1+\cdots+ m_lV_l$ be generic. Then, as it is well known (see, e.g., \cite{kac}), ${\rm Ext}(V_i,V_j)=0$ for $1\leq i\neq j \leq l$. It is also well known that for any representations $U,V$ holds:
\begin{equation}\label{homext}
\langle \dim U, \dim V\rangle = \dim {\rm Hom}(U,V) - \dim {\rm Ext}(U,V).
\end{equation}  
In particular, we have $\langle \beta_i,\beta_j\rangle=\dim{\rm Hom}(V_i,V_j)$. However, by our assumption, $X$ is locally semi-simple and, similar to  simple representations, homomorphisms between non-isomorphic locally simple factors vanish, see \cite[Proposition~8]{sh07}. Therefore, by formula (\ref{homext}), $\langle \beta_i,\beta_j\rangle=0$. Conversely, if the equation $\langle \beta_i,\beta_j\rangle=0$ for every pair of different summands of generic decomposition, this decomposition is also locally semi-simple as follows from \cite[Theorem~4.3]{sh09}. Indeed, in terms of that Theorem, almost-loopless condition ($m_i=1$ if $\langle\beta_i,\beta_i\rangle <0$) follows from the properties of generic decompositions, the {\it local quiver} has no edges except for loops and the number of summands is the same as for generic (as it is generic itself). 

So far we proved 
${\bf c}\Leftrightarrow {\bf a}\Rightarrow{\bf b}\Leftrightarrow{\bf d}$. We conclude the proof of Theorem with the implication ${\bf d}\Rightarrow{\bf c}$. To see that, we recall that the equation $p\beta_j[0]=n(\beta_j[1]+\cdots+\beta_j[k])$ means that the so-called {\it Schofield weight} $\sigma_{\alpha}=[-p,n,\cdots,n]$ vanishes on every summands $\beta_j$ of the generic decomposition for $\alpha$ as it obviously vanishes on $\alpha$. Schofield weight has a universal definition as $\sigma_{\alpha}(\beta)=\langle \alpha, \beta\rangle - \langle \beta, \alpha\rangle$ (see \cite[Definition~3.9]{fm}) and this definition also makes $\sigma_{\alpha}(\alpha)=0$ obvious. Now, assuming $\langle \beta_i, \beta_j\rangle=0$ as in {\bf d}, we conclude:

$$
\langle \alpha, \beta_i\rangle=\sum_{j=1}^l \langle m_j\beta_j, \beta_i\rangle=m_i\langle \beta_i, \beta_i\rangle=\sum_{j=1}^l \langle \beta_i,m_j\beta_j\rangle= \langle  \beta_i,\alpha\rangle \qed
$$

\section{Flip-Flop Algorithm for iPCA}\label{FFS}

As we mentioned above, the authors of \cite{ta} gave a simple Flip-Flop algorithm for (unpenalized MLE of) iPCA but didn't believe it works in any meaningful generality. The input is $X=(X_1,X_2,\cdots,X_k)$ and the output are estimators for covariance, symmetric (or hermitian) $n\times n$ matrix $\Sigma$ and $p_j\times p_j$ matrices $\Delta_j, j=1,\cdots,k$. The algorithm is so simple that we find it important to give it:

{\bf 1:} Initialize $\Sigma,\Delta_j,j=1,\cdots,k$ with identity matrices

{\bf 2: while} not converge, do:

{\bf 3:} \hspace{0.2cm} Update $\Sigma = \frac{1}{p}\sum_{j=1}^k X_j \Delta_j^{-1}\overline{X_j^{\top}}$

{\bf 4:} \hspace{0.2cm} {\bf for} $j=1,\cdots,k$

{\bf 5:} \hspace{0.4cm} Update $\Delta_j=\frac{1}{n} \overline{X_j^{\top}} \Sigma^{-1} X_j$

We  implemented this algorithm in \cite{ff} and did experiments with random $X$ and different dimension vectors $\alpha$. Our experiment fully confirmed the statement of Theorem \ref{main_result}. Namely, once $\alpha$ meets the conditions of Theorem \ref{main_result} and we select a generic $X\in R(Q_k,\alpha)$ at random, the Theorem guarantees that $X$ is likely to be locally semi-simple, and we see in our experiments that  algorithm converges. Conversely, if this is not the case for $\alpha$, the algorithm diverges and the Frobenius distance between the points before and after an iteration increases fast. Moreover, talking in more software language, in such bad cases,  at some step,  after an update in line {\bf 3}, $\det(\Sigma)$ becomes too small, so the inversion in the next call of line {\bf 5} throws an exception visible for a user.

\section{Conclusions and next steps}

For an iPCA setting $\alpha=[n,p_1,\cdots,p_k]$ we gave a criterion, Theorem \ref{main_result} for MLE to exist for generic $n$ observed vectors with feature groups of sizes $p_1,\cdots,p_k$. For practical usage of this criterion we propose our two implementations of Derksen-Weyman algorithm, \cite{dw}. First, a  software with GUI \cite{t} that is easy in use and may be interesting for, first of all, quiver theorists, because it gives more than just Derksen-Weyman algorithm. Second, we shared C++ source code needed to  run Derksen-Weyman algorithm in \cite{ff}. Also \cite{ff} contains an implementation of Flip-Flop algorithm and a demo test function that iterates over  dimension vectors $\alpha$ for a star quiver, selects those fitting to the conditions of \ref{main_result}, and for such $\alpha$, samples data at random, applies Flip-Flop algorithm and verifies convergence. 

Therefore, MLE exists for much more cases of iPCA than \cite{ta} believed and the next step is to use it in practical experiments with big data. However, for big data the convergence speed of Flip-Flop is too slow and also it uses matrix inversion that is time consuming. So we address a problem of finding fast MLE algorithms for practical needs.

\end{document}